\documentclass{siamart190516}
\usepackage{amsmath,amsfonts,color,ifpdf}
\usepackage{algorithm,algorithmic}
\usepackage{epsfig, graphicx}
\usepackage{cite}

\newcommand{\btrue}{b_{\mathrm{true}}}
\newcommand{\xtrue}{x_{\mathrm{true}}}

\title{Multigrid preconditioning for regularized least-squares problems\thanks{Submitted to the editors DATE.\funding{Kilmer’s work on this project has been partially supported by NSF DMS-1720291, NSF DMS-1821148, and NSF HDR grant CCF-1934553.  MacLachlan's work has been partially supported by an NSERC discovery grant.}}}
\author{Matthias Bolten  \and Scott P. MacLachlan \and Misha E. Kilmer}

\begin{document}
\maketitle

\begin{abstract}
In this paper, we are concerned with 
efficiently solving the sequences of regularized linear least squares problems associated with employing
Tikhonov-type regularization with regularization operators designed to enforce edge recovery.  
An optimal regularization parameter, which balances the fidelity to the data with the edge-enforcing constraint term, is 
typically not known a priori.  This adds to the total number of regularized linear least squares problems that must be solved before the final image can be recovered.  Therefore, in this paper, we determine effective multigrid preconditioners for these sequences of systems.  We focus our approach on the sequences that arise as a result of the edge-preserving method introduced in \cite{Gazzola_etal_2020}, where we can exploit an interpretation of the regularization term as a diffusion operator; however, our methods are also applicable in other edge-preserving settings, such as iteratively reweighted least squares problems. Particular attention is
paid to the selection of components of the multigrid preconditioner in order to achieve robustness for
different ranges of the regularization parameter value.    In addition, we present a parameter culling approach that, when used with the L-curve heuristic, reduces the total number of solves required.  We demonstrate 
our preconditioning and parameter culling routines on examples in computed tomography and image deblurring.

%
%
\end{abstract}



\section{Introduction}  
In discrete linear inverse problems such as image reconstruction and image restoration, 
the forward operators are known to be ill-conditioned.  
Due to the presence of noise in the data, regularization is required to produce a meaningful
image from the measured data.   
In this paper, we are interested in designing efficient solvers for the sequence of regularized least squares problems
\begin{equation} \label{eq:LSreg}  \min_{x} \| A x - b \|_2^2 + \lambda^2 \| M^{(\ell)} x \|_2^2 \end{equation}
or, equivalently,
\begin{equation} \label{eq:LStall} \min_x \left\| \left[ \begin{matrix} A \\ \lambda M^{(\ell)} \end{matrix} \right] x - \left[ \begin{matrix} b \\ 0 \end{matrix} \right] \right\|_2 ,\end{equation}
where $b$ represents the (noisy) data, $x$ is the recovered solution (e.g., pixel-wise values of an image), $A$ is the so-called ``forward operator'', mapping an image to its data, and the matrices $M^{(\ell)}$ denote regularization operators.  The scalar $\lambda$ is called the {\it regularization parameter}, and serves to balance the data fidelity term with the regularization.  We are particularly interested in regularization operators designed to 
recover edge information in restored or reconstructed images.  
Routinely, this class of regularization operators are expressed as 
\begin{equation} \label{eq:ml} M^{(\ell)} = D^{(\ell)} L,  \end{equation}
where $L$ is the discrete gradient matrix and $D^{(\ell)}$ denotes a
diagonal weighting matrix that is updated in an iterative process as additional information about the edges becomes available, as denoted by index $\ell$.   
We refer to this iteration over $\ell$ as the {\it outer iteration}.  

Such sequences of problems
arise frequently in situations where the entries in the diagonal weighting matrix would ideally be computed from the desired $x_{true}$ but, as this is unknown, the weights are recomputed adaptively from increasingly better estimates of the solution. 
One such approach, which is the method of focus in the present paper,
is described in \cite{Gazzola_etal_2020}.  In that paper, the diagonal
weighting matrices, $D^{(\ell)}$, are updated from iteration to
iteration to preserve edges by reducing values in $D^{(\ell)}$ that
correspond to edges detected in that iteration,  
such that each new weighting matrix is a diagonal weighting of the previous weighting matrix.   However, other examples
of weighting matrices with properties similar to those we exploit here abound in the literature, most notably in iterative quadratic approximations to the total variation regularization operator or to the $p$-norm of the gradient (see, for example, \cite{IRNekki,renaut2017hybrid,IRNtv,wohlberg2008lp}); therefore, we expect our method to be more widely applicable. 
%


There are many heuristics in the literature for choosing the optimal
regularization parameter, $\lambda^{*,\ell}$ (see \cite{hansen2010discrete} and references therein).
While {our solvers would be applicable in scenarios where other selection mechanisms are used}, 
for ease of presentation, we will consider only
one well-known heuristic for estimating $\lambda^{*,\ell}$ for fixed $\ell$,
called the L-curve criterion.   The main consideration is that most
heuristics that compute an estimate of the optimal parameter,
including the L-curve, require measurements of the data-fidelity and or
constraint terms as a function of the parameter, $\lambda$.  Thus,
even when done carefully, much of the computation involved in solving
the inverse problem goes into the effort of multiple solves for each
outer iteration, $\ell$, as a function of $\lambda$. 
The L-curve \cite{HansenOleary93} is a parametric plot of $( \log_{10}
\|A x_{\lambda} - b \|_2,\log_{10} \| M^{(\ell)} x_{\lambda} \|_2)$ as
a function of $\lambda$, where $x_{\lambda}$ denotes the solution of~\eqref{eq:LSreg} at outer iteration $\ell$ and parameter $\lambda$.   The general idea is that the point of
maximum positive curvature in this curve occurs where there is a good balance between the effect of the regularization term and that of the data-fidelity term.   If several discrete values of $\lambda$ are given, then the discrete value that lies closest to the corner of the interpolant is taken as the estimate of the optimal value.  

The naive approach to generating points on the L-curve sufficient to make a determination about the value in the corner would be to solve (\ref{eq:LSreg}) (or its equivalent form in~\eqref{eq:LStall}) for a few choices of $\lambda$ and then plot the curve. 
Due to the size and structure of the operators, iterative Krylov subspace solvers such as CGLS or LSQR \cite{PaSa82a,PaSa82b} are usually the 
methods of choice for solving~\eqref{eq:LSreg} for each data point. 
 The regularized solution at the $\ell$th outer iteration would correspond to the solution at the maximum positive curvature point on the curve.   If this approach is implemented as stated, however, it is expensive for the following reasons:
   \begin{itemize}
      \item When $\lambda$ is small,   $\left[ \begin{matrix} A \\ \lambda M^{(\ell)} \end{matrix} \right] $ is likely to be ill-conditioned, as the emphasis is on the ill-conditioned top block.   Thus, many Krylov iterations are required for an {\it accurate} solution.  Accuracy is key, since inaccurate solutions deteriorate the shape of the L-curve and allow miscalculations of the corner.
      \item When $\lambda$ is large, the emphasis is on the bottom block, but as $\ell$ increases, $M^{(\ell)}$ tends
      to also be ill-conditioned due to more accurate diagonal weighting and near zero terms.
      \item A new L-curve must be determined for each outer index $\ell$ since the optimal choice of $\lambda$ changes with $\ell$.
      \end{itemize}
      
As an alternative to the naive approach outlined above, hybrid iterative methods for the class of regularization operators we consider (namely, those having $M^{(\ell)}$ that is rectangular with a possible non-trivial null-space) have been proposed in the literature.  In such methods, the regularized approximation to each $x_\lambda$ is constrained to live in the same $k$-dimensional subspace whose basis is generated by an {\it inner iterative procedure}  \cite{Kilmer_Hanson_Espanol2007,Gazzola_etal_2020}.  Related methods for the Tikhonov regularized problem, which are not immediately applicable for the case of rectangular $M^{(\ell)}$ we consider here, are given in, for example, \cite{GaNa14,GaNoRu15,KiOLe01,OlSi81,ChPa15}.   Hybrid
methods have an advantage that, given a suitable $k$, an appropriate choice for $\lambda$ can be chosen for the $k$-dimensional projected problem very cheaply, and only one subspace is created for all values of $\lambda$.   

The cost of the hybrid methods varies by algorithm, but the main expense comes from computing the subspace into which the problem is projected, which in turn requires a minimum of $k_\ell$ matrix-vector products, where $k_{\ell}$ is the dimension of the subspace at outer iteration $\ell$.     
  The hybrid method \cite{Kilmer_Hanson_Espanol2007} that is
used in \cite{Gazzola_etal_2020}, though very effective at jointly predicting the regularization parameter for the $\ell$th system and producing the corresponding regularized solution, can require many more than $O(k_{\ell})$ matrix-vector
products in producing the $k_\ell$ dimensional subspace, for the following reason.  In order to generate the basis for the subspace, there are repeated inner loop calls to LSQR (itself an iterative method requiring matrix-vector products) which are required to compute approximate orthogonal projections at each step.
In other words, each iteration of the hybrid method makes two calls to LSQR, and LSQR requires some number of matrix-vector products.  Thus, one step of the hybrid regularization method costs a multiple of the cost of a matrix-vector product, and that multiple may be (much) bigger than $O(1)$.  


Our goal in this paper is to propose an alternative that more closely follows the outline of the so-called naive approach, but which is computationally efficient in computing the desired regularization solution.  We apply our method to the sequence of regularized systems with regularization operators $M^{(\ell)}$ produced by the approach in \cite{Gazzola_etal_2020}.   As noted above, in this case $M^{(\ell)} = D^{(\ell)} L$ where $L$ is the discrete gradient and $D^{(\ell)}$ is a diagonal weighting matrix that is generated from gradient information from previous regularized solutions and designed in such a way as to preserve edges while removing noise.  
 Our goal is to gain computational efficiency by reducing the cost of the solves {\it as well as} reducing the total number of overall solves required to find the final regularized solution.   
To address the solver cost, we exploit the fact that the regularization term as it appears in the regularized normal equations corresponding to (\ref{eq:LStall}) resembles a 
diffusion-like operator (see also \cite{CRVogel_2002a} for similar insight into such regularization terms).  
To address the total number of solves required, we leverage the observation from \cite{Gazzola_etal_2020} that the optimal regularization parameters increase, up to a point, as $\ell$ increases (i.e. as additional edge information becomes visible in the regularized solutions).     

Our specific contributions are:
   \begin{enumerate}
     \item FGMRES-accelerated AMG for solving the system using the normal equation with few iterations when $\lambda$ is larger and favors the diffusion-like term.
     \item Non-exact solves with $\lambda$ small, just to know enough about the L-curve to find the corner.
     \item Leveraging properties of the L-curves as functions of $\ell$ to prune the L-curves  and isolate the best choice of
$\lambda$ for a given $\ell$.
     \item A new heuristic for stopping the outer iteration on $\ell$ to ensure that our algorithm is automatic and free from the need for user-based tuning. 
     \end{enumerate} 

This paper is organized as follows.   The background and motivation is discussed in more detail in \Cref{sec:problem}.  In \Cref{sec:firstimprove}, we detail our approach: the multigrid solver is covered in \Cref{subsec:amg_for_the_regularized_linear_systems} with the pruning routine for the outer iterations in  \Cref{sec:secondimprove}.
The efficiency of our proposed method is demonstrated in Section \Cref{sec:numerical}.  Finally, \Cref{sec:conclusions} contains conclusions and future work.

\section{Background and Problem Formulation} \label{sec:problem}
We consider the solution of discrete ill-posed problems with linear
forward models of the form
\begin{equation}
\label{eq:exact}
A x = b = \btrue + \eta,
\end{equation}
where $A \in \mathbb{R}^{m\times n}$ 
is the forward
operator, $\xtrue$ is the ``true'' solution that we look to recover,
$\btrue = A \xtrue$ is the exact, noise-free data generated by applying $A$
to $\xtrue$, and $\eta$ is an unknown additive white noise vector.  We
assume that both the true solution and noise-free data are unknown, so
that we can only access the problem via the forward operator, $A$, and
the noisy data vector, $b$.  We focus on imaging problems, where $x$ and
$\xtrue$ are taken to be the vectorized forms of digital images, and
$A$ has singular values that decay rapidly towards zero.  As is typical with discrete ill-posed problems (and in contrast with
discretized PDE operators) the
singular vectors corresponding to the largest singular values
represent smooth modes and those corresponding to the smallest
singular values represent high-frequency modes.  If $\text{rank}(A)< n$, then, independent of the presence of the noise, there is no unique least-squares solution.  However, even when $\text{rank}(A) = n$, the (unique)
least squares solution is problematic:
due to the
presence of the noise and the decay of the singular values, the (least norm) linear
least squares solution 
will be contaminated by noise, with the extent of noise corruption
determined by the decay rate of the singular values as well as the
noise level \cite{Hansenbk}. 


\typeout{here2}

In applications where image edge-preservation is required, the regularized problem should encourage 
solutions that have edges preserved. 
Thus, the least squares problem is replaced by a regularized 
least squares problem in order to simultaneously combat the (near) rank deficiency and
noise: 
\[ \min_x \| A x - b \|_2^2 + \lambda^2 R(x), \]   
where the regularization functional should be designed to encourage solutions with edges.  

Notationally, let $x$ denote the vectorized form of an $n_v \times n_h$ image.
The matrix
\[ L_v = \begin{bmatrix} -1 & 1 & 0 & 0 & \cdots \\ 0 & -1 & 1 & 0 & \cdots \\
                                  0 & 0 & -1 & 1 &  \cdots \\
                                  \vdots & \ddots & \ddots & \ddots &\ddots \\
                                  0 & 0 & \cdots & -1 & 1 \end{bmatrix} \in \mathbb{R}^{(n_v-1) \times n_v},  \]
represents the discrete first derivative operator acting on a vector
with $n_v$ entries.   Likewise, we can define $L_h \in
\mathbb{R}^{(n_h-1) \times n_h}$.   The discrete gradient operator is, then, given by
\begin{equation}\label{eq:gradient}
L = \begin{bmatrix} I \otimes L_v \\ L_h \otimes I \end{bmatrix} \,.
\end{equation}
It is a straightforward exercise to show that 
$Lx \in \mathbb{R}^{n_h(n_v-1)+(n_h-1)n_v}$ contains the edge information of $x$, meaning that entries of $Lx$ are large in positions corresponding to edges in the image, where the values of $x$ change between adjacent pixels.

 It is well known that using $R(x) = \| x \|_2^2$ or $R(x) = \| L x \|_2^2$ will smooth edges, rather than enhance them.  On the other hand, these choices of regularization have the advantage that we need to only solve a single quadratic problem (assuming $\lambda$ is known).  
For edge-preservation, typical choices of $R(x)$ are (isotropic or non-isotropic) total variation of $Lx$, or $\| L x \|_p^p$, where $p$ is close to 1; with these, solving the Tikhonov edge-regularized problem requires solving a sequence of quadratic problems.
Consider, for example,  the \textit{iteratively reweighted norm} (IRN) method, proposed
in \cite{IRNekki} and further developed in
\cite{wohlberg2008lp,IRNtv,renaut2017hybrid}. In this approach, a non-quadratic
constraint, such as $\|Lx\|_p^p$ for $p\approx 1$ or the total
variation of $x$, $\mathrm{TV}(x)$, is reformulated as a sequence of
quadratic regularization problems,
\begin{equation} \label{eq:opt}
\min_{x} \| A x - b \|_2^2 + \lambda^2 \| M^{(\ell)} x \|_2^2,
\end{equation}
where $M^{(\ell)}$ is defined based on the discrete gradient matrix ,$L$, and the approximate solution
of the $(\ell-1)$th problem. 

In the present paper, we instead use the sequence $M^{(\ell)}$ of regularization operators 
designed for
edge preservation introduced in \cite{Gazzola_etal_2020}.  
{\it In the description of their weighting method below, 
we have assumed that near optimal values $\lambda^{*,\ell}$ are known}, so that we have decoupled the description of the weighting process from the parameter selection process.   However, a key 
issue in the computational complexity associated with finding the ultimate regularized solution is that these values
are, in fact, {\it not known a priori}.   The parameter selection costs and mitigation approach we present are described in the next section. 

 As noted in the introduction, the regularization 
operator from \cite{Gazzola_etal_2020} takes the form $M^{(\ell)} = D^{(\ell)} L$
where $L$ is as defined above and $D^{(\ell)}$ is a diagonal weighting matrix that is updated as follows. 
First, let $D^{(0)} = I$, and denote by $x^{(*,0)}$ the solution to 
\[ \min_x \| A x - b \|_2^2 + (\lambda^{*,0})^2 \| D^{(0)} L x \|_2^2, \]
where $\lambda^{*,0}$ is the estimated optimal regularization parameter for this problem.  
We next set $g^{(1)} = | L x^{(*,0)} | ./ \| L x^{(*,0)} \|_\infty$, and 
$d^{(1)} := {\bf 1} - (g^{(1)}).^q$, for $q > 0$ fixed at the same value for all iterations, $\ell$.    
Here, the absolute value is meant to denote an element-wise operation on the vector argument.  Likewise, the ``dot" in front of the division and exponentiation characters is meant to denote a componentwise operation, and ${\bf 1}$ is the vector of all ones.
%
Now, let $D^{(1)} = \text{diag}(d^{(1)}) D^{(0)}$, and
the next system in the sequence to be solved is 
\[ \min_x \| A x - b \|_2^2 + (\lambda^{*,1})^2 \| D^{(1)} L x \|_2^2. \]
The next system is generated from the solution to this problem by defining $g^{(2)} := | D^{(1)} L x^{(*,1)} | ./ \| D^{(1)} L x^{(*,1)} \|_\infty$, then setting $d^{(2)} = {\bf 1 }  - (g^{(2)}).^q$, and by updating the diagonal matrix $D^{(2)} = \text{diag}(d^{(2)}) D^{(1)}$.  
In general, the diagonal matrix is updated as $D^{(\ell)} =
\text{diag}(d^{(\ell)})D^{(\ell-1)}$, 
where we are defining
\begin{equation} \label{eq:g}
g^{(\ell)} = |D^{(\ell-1)} Lx^{(*,\ell-1)}| . / \| D^{(\ell-1)} Lx^{(*,\ell-1)}\|_\infty, \qquad d^{(\ell)} = {\bf 1} - (g^{(\ell)}).^q . 
\end{equation}



In \cite{Gazzola_etal_2020}, they employ 
the iterative, joint bidiagonalization technique from
\cite{Kilmer_Hanson_Espanol2007}. 
Specifically, regularized solutions for any $\lambda$ are obtained by solving
\begin{equation} \label{eq:restricted} \min_{x \in \mathcal{Z}_k } \| A x - b \|_2^2 + \lambda^2 \| M^{(\ell)} x \|_2^2, \end{equation}
where $\mathcal{Z}_k$ is a $k$-dimensional subspace generated iteratively. 
The idea of this technique is to create a basis set,
$\{z^{(\ell)}_j\}_{j=1}^k$, such that
\[
AZ_k^{(\ell)} = U_{k+1}^{(\ell)}B_k^{(\ell)} \text{ and
}M^{(\ell)}Z_k^{(\ell)} = \hat{U}_k^{(\ell)}\hat{B}_k^{(\ell)},
\]
where $Z_k^{(\ell)}$ is the matrix with columns given by
$z^{(\ell)}_j$, $U_{k+1}^{(\ell)}$ and $\hat{U}_k^{(\ell)}$ are
orthogonal matrices, and $B_k^{(\ell)}$ and $\hat{B}_k^{(\ell)}$ are
bidiagonal matrices.  Note that this basis {\it does not depend on $\lambda$ but does depend on $M^{(\ell)}$}. 
Since the regularized solutions are restricted to lie in
the range of $Z_k^{(\ell)}$, we have that
$x_{\lambda} = Z_k^{(\ell)} w_{\lambda}$, which leads to the
$k$-dimensional projected problem
\begin{equation}
  \label{eq:sequence_hybrid}
\min_{w_\lambda} \| B_k^{(\ell)}w - \beta_1 e_1 \|_2^2 + \lambda^2 \| \hat{B}_k^{(\ell)} w \|_2^2,
\end{equation}
where $\beta_1 e_1 = \left(U_{k+1}^{(\ell)}\right)^Tb$.   

Here, $k$ is small relative to the problem size.   
Now $\lambda^{*,\ell}  \in \{ \lambda^{(i)} \}_{i=1}^m$ is estimated to be the optimal regularization parameter for the small projected problem, and the solution is updated as $x^{(i,\ell)} = Z_k^{(\ell)} w^{(i,\ell)}$ where $w^{(i,\ell)}$ is the solution to (\ref{eq:sequence_hybrid}) for
$\lambda^{*,\ell} = \lambda^{(i)}$.   
Finding the L-curve and its corner for the size-$k$ projected problem { is very cheap}, because there are recurrence relations for the residual and constraint terms as functions of $\ell$.  Likewise, the $x^{(i,\ell)}:=Z_{k}^{(\ell)} w^{(i,\ell)}$ can be computed each iteration $k$ for all of the $m$ values of $\lambda^{(i)}$ simultaneously using short term recurrences for each solution\footnote{If storage is at a premium, an alternative is to run the iteration once to find $\lambda^{*,\ell}$, then run it again with the selected parameter to find the corresponding solution.}. To be clear, there is a different L-curve for each iteration index $k$.  However, once $k$ is sufficiently large, the curves and, particularly the corner values, will stagnate, indicating that
there is no use in iterating further. We call the final iteration value $k_\ell$, indicating that the final required number of iterations needed depends both on properties of $A$ and the regularization operator and therefore can take on different values as functions of $\ell$. 

On the face of it, this seems a very efficient approach in terms of computations and storage to solving the desired problem: for a fixed outer iteration $\ell$, generate only one subspace of dimension $k_\ell$; find suitable $\lambda^{*,\ell}$ for the projected problem and corresponding $x^{(*,\ell)}$ from recurrence relations.  In comparison, the naive approach would require running an iterative method such as LSQR for each shifted system (i.e. $m$ different subspaces of varying dimension) to get the data to form an L-curve.  
However, as we discuss in more detail in the next section, {\it there are some hidden, non-trivial costs involved in the hybrid approach that must be considered.}    
It is the overall expense of both extremes (i.e. the hybrid and naive approaches) that is driving us to consider the alternative algorithm that we develop in this paper. 


%

\section{Improving Efficiency of Inner-Outer Iterations} \label{sec:firstimprove}

Recall that $\ell$ refers to the outer-iteration and denotes a change in the regularization operator.   Thus far, we have discussed two choices for solving the linear system at each outer iteration:
\begin{enumerate}
\item Use a hybrid approach (where all solutions to linear systems lie in the same subspace) to find the regularization parameter and corresponding solution, or\item Use the naive approach of independently solving all $m$ systems with an iterative method (e.g. LSQR) for each regularization parameter to tolerance (noting that each solution now lies in a different subspace), form the L-curve, select the parameter and corresponding solution.
\end{enumerate}
We first compare these approaches in terms of the cost of matrix-vector products with the stacked matrix.  

For the hybrid approach, even when the dimension, $k_\ell$, of the subspace for the 
regularized, projected problem is relatively small, the algorithm may require {\it far more than $k_\ell$ matrix-vector products}. This is because each iteration of the hybrid algorithm requires two ``interior'' calls to LSQR, in order to generate the approximate orthogonal projections needed to compute the columns of the basis matrices.   Depending on the properties of the operators $A$ and $M^{(\ell)}$, LSQR may require many iterations (say, $j_k$) to compute the desired projection.  Thus, the cost of $k_\ell$ steps of the hybrid algorithm to select the best regularization parameter (on the projected problem) and compute the corresponding regularized solution is dominated by $O(\sum_{k=1}^{k_\ell}  j_k)$ matrix-vector products.    Typically, there is not much variability in $j_k$ across the hybrid iterations, so the dominant cost to execute the hybrid algorithm can be approximated by that of $O( k_\ell J)$ matrix-vector products.  
   
On the other hand, the hybrid method has the advantage of very efficiently choosing the regularization parameters. 
We can compare this to the naive approach of using an iterative method to compute (\ref{eq:LStall}) for each new value of $\lambda^{(i)}, i=1,\ldots,m$, {\it to a tolerance sufficient to be able to differentiate the effect of the regularization parameter}, so that we can employ a heuristic such as the L-curve to select a parameter and get our solution.  
Unfortunately, the number of iterations required to solve the $i$th system to tolerance may strongly depend 
on $\lambda^{(i)}$, as this value affects the conditioning.  
The dominant cost of an iteration of LSQR is the cost of the matrix-vector product.  Thus, if
system $i$ requires $\kappa_i$ iterations to reach convergence, 
then the dominant cost to be able to select our parameter and generate the corresponding solution 
is $O(\sum_{i=1}^m \kappa_i)$ matrix-vector products.  
For the smallest shifts, the system is typically very ill-conditioned, and $\kappa_i$ can be very large.  Furthermore, as the outer iterations progress, the systems for the largest shifts will be ill-conditioned as well.    
Thus, the naive approach can be much more expensive than than the hybrid method, particularly when considered over all the outer iterations.   
  
Therefore, we address the two issues that keep the naive approach from being competitive with the hybrid approach.
To address
the first issue of the number of iterations, we introduce a multigrid preconditioner for the solution
of the problems.  To address the second, we
introduce two key improvements in the selection of $\lambda^{*,\ell}$, by
solving over fewer potential values for the regularization parameter
and by using the solution to the system (at the current $\ell$) for the neighboring $\lambda$ as an  
initial guess for the next system solve. 


\subsection{Multigrid}\label{subsec:multigrid}
Multigrid methods (see, for
example, \cite{WLBriggs_VEHenson_SFMcCormick_2000a, UTrottenberg_etal_2001a}) are optimal solvers for systems arising from the
discretization of elliptic PDEs. Geometric multigrid methods make use
of the following observation: When an iterative method like (weighted)
Jacobi or Gauss-Seidel is applied to a linear system
\[
  K x = b, K \in \mathbb{R}^{n \times n}, x, b \in \mathbb{R}^n,
\]
arising from the discretization of a simple elliptic PDE, the
reduction of the error is relatively slow and depends on the
discretization parameter, $h$. Nevertheless, the error $e = K^{-1} b -
x$ is smooth after only a few iterations of the iterative
method. Obviously, a smooth error can be represented well on a coarser
mesh. Using the definition of the residual $r = b - K x$, we can obtain
the error as the solution of $K e = r$. If coarse representations,
$K_c$ of $K$ and $r_c$ of $r$, are available, then the system $K_c e_c = r_c$
can be solved instead and the solution is less expensive. The
restriction of $r$ to $r_c$ can be obtained by simple injection or by
more complicated schemes such as, e.g., full weighting. The coarse-grid
operator, $K_c$, can be obtained using a rediscretization of the
continuous operator on the coarser mesh, or by using a variational
principle resulting in the Galerkin operator. The error, $e_c$, on the
coarse mesh is then interpolated to the fine one, yielding an
approximation $\tilde{e}$. Finally, the current approximate solution
can be updated as $x = x + \tilde{e}$. As the interpolation process
introduces high-frequency error components in the approximate
solution, in most cases a couple of additional relaxation steps are
applied. The resulting method is called a
\emph{two-grid} method.

This idea can be applied recursively on each level, with level
$k$ approximating the solution of the coarser system $K_{k+1} e_{k+1} = r_{k+1}$,
using $\gamma$ iterations of the described method on
level $k+1$, resulting in a multigrid method. Given relaxation
schemes on level $k$ specified by update operators,
$\mathcal{S}_k$ and $\tilde{\mathcal{S}}_k$,
$k = 0,\dots,k_\text{max}$, number of pre- and post-relaxation steps, $\nu_1$ and $\nu_2$, and restriction and prolongation operators,
$R_k$ and $P_k$, $k = 1,\dots,k_{\text{max}-1}$, a generic
multigrid algorithm can be found in Algorithm~\ref{alg:mg}.

\begin{algorithm}
  \caption{Multigrid cycle $x_k = \mathcal{MG}_k(x_k,b_k)$}
  \label{alg:mg}
  \begin{algorithmic}
    \STATE{$x_k \leftarrow
      \mathcal{S}_k^{\nu_1}(x_k,b_k)$}
    \STATE{$r_k \leftarrow b_k - K_k
      x_k$}
    \STATE{$r_{k+1} \leftarrow R_k r_k$}
    \STATE{$e_{k+1} \leftarrow 0$}
    \IF{$k+1 = k_\text{max}$}
    \STATE{$e_{k_\text{max}} \leftarrow
      K_{k_\text{max}}^{-1} r_{k_\text{max}}$}
  \ELSE
  	\FOR{$j = 1,\dots,\gamma$}
  		\STATE{$e_{k+1} \leftarrow \mathcal{MG}_{k+1}(e_{k+1},r_{k+1})$}
	\ENDFOR
  \ENDIF
  \STATE{$e_k \leftarrow P_k e_{k+1}$}
  \STATE{$x_k \leftarrow x_k +
    e_k$}
  \STATE{$x_k \leftarrow
    \tilde{\mathcal{S}}_k^{\nu_2}(x_k,b_k)$}
  \end{algorithmic}
\end{algorithm}

When the coarsening is not easily determined geometrically, an
algebraic multigrid method (AMG) can be
used \cite{ABrandt_SFMcCormick_JWRuge_1984a, JWRuge_KStuben_1987a,
KStuben_2001a}. For a given relaxation scheme, classical AMG
automatically selects coarse mesh points and the necessary restriction
and prolongation operators, using heuristics that are
typically derived based on M-matrix or similar properties. AMG for
symmetric problems usually chooses $R_k = P_k^T$, and the coarse-grid
operator is chosen to be the Galerkin coarse-grid operator given by
$K_{k+1} = P_k^T K_k P_k$.  Once this ``setup phase'' is complete, the
AMG solution phase proceeds exactly as in the geometric MG setting, as
described in Algorithm \ref{alg:mg}.

In what follows, we make use of components of the classical
(Ruge-St\"uben) AMG algorithm as described
in \cite{JWRuge_KStuben_1987a} and elsewhere, but within a somewhat
nonstandard setup and solve phase, which are adapted to the problem
setting at hand.  We first describe the classical components.  Given a
matrix, $K$, the coarse grid is selected
using a two-pass algorithm that partitions the degrees of freedom in
$K$ into two disjoint sets, denoted as $C$ and $F$.  First, a maximal
independent set of the graph of the so-called strong connections in
matrix $K$ is formed.  This amounts to first filtering the matrix $K$
to remove entries where
\[
-k_{ij} \leq \theta \max_{m\neq i} \{-k_{im}\},
\]
where we take $\theta = 0.25$, as is commonly done.  The maximal
independent set is then formed by a greedy algorithm that sequentially
identifies the unsorted point with largest measure (initialized as
number of strong connections), marks this as a $C$-point, then marks its
as-yet unsorted strong neighbors as $F$-points and increments the
measure of the unsorted strong neighbors of the $F$-points to make
these more attractive for subsequent selection as $C$-points.  Given
this tentative partitioning, a second pass is performed, where
additional points are moved from the $F$-set to the $C$-set, in order to
satisfy heuristic measures that indicate an effective interpolation
operator can be found.  Once the second pass has finalized the
partitioning, the interpolation operator, $P$, is formed.  For each
point in $C$, we assign a coarse-grid index to the point, and
interpolate directly from the coarse-grid index to its corresponding
point in $C$.  For interpolation to points in $F$, we make the assumption
that we are interpolating errors that vary smoothly between strongly
connected points, leading to an interpolation formula for each $i\in
F$ of $e_i = \sum_{j\in C_i} w_{ij}e_j$, where $C_i$ is the set of
points in $C$ that are strongly connected to $i$.  The interpolation
weights are given by
\[
w_{ij} = - \frac{k_{ij} + \displaystyle\sum_{m \in
    F_i^s}\frac{k_{im}k_{mj}}{\sum_{l \in C_i} k_{ml}}}{k_{ii} +
    \displaystyle\sum_{m \in F_i^w} k_{im}},
\]
where the set $F_i^s$ is the set of points in $F$ that are strongly
connected to $i$, while the set $F_i^w$ is defined so that all points
adjacent to $i$ are contained in the disjoint union $C_i \cup
F_i^s \cup F_i^w$.

\subsection{AMG for the regularized linear systems}\label{subsec:amg_for_the_regularized_linear_systems}
To clarify notation, we recall that $x^{(i,\ell)}$ is the vector that minimizes \eqref{eq:LStall}.
The normal equations are readily given by
\begin{equation}
\label{eq:normal}
\left(A^TA + \left(\lambda^{(i)}\right)^2 \left(M^{(\ell)}\right)^TM^{(\ell)}\right)x^{(i,\ell)} = A^Tb.
\end{equation}
As mentioned above, the first matrix in \eqref{eq:normal} corresponds to the ``data'' term
from the inverse problem, which is often an integral operator or
non-differential operator that is difficult to treat directly using
AMG due to its spectral properties and expected large number of
nonzero entries per row.  In contrast, the second matrix
in \eqref{eq:normal},
\[
\left(M^{(\ell)}\right)^TM^{(\ell)} = L^T\left(D^{(\ell)}\right)^2L,
\]
represents a variable-coefficient diffusion operator, for which AMG
techniques are expected to provide quite effective and efficient
solvers.

As a first observation, we note that for large-enough values of
$\lambda^{(i)}$, the system will be dominated by its second component.
This motivates the design of an AMG-based algorithm to
solve \eqref{eq:normal}.  However, directly applying AMG to $A^TA
+ \left(\lambda^{(i)}\right)^2 L^T\left(D^{(\ell)}\right)^2L$ is
unlikely to be effective, due to the cost of applying the underlying
graph algorithms to the much denser matrix coming from $A^TA$.  Here,
we propose an alternative approach, that takes advantage of the better
sparsity of the rectangular operator in \eqref{eq:LStall} over the
normal operator in \eqref{eq:normal}.  To define the coarse-grid
problem, we first form the second operator in \eqref{eq:normal},
$L^T\left(D^{(\ell)}\right)^2L$, and apply the classical AMG algorithm
discussed above (and in \cite{JWRuge_KStuben_1987a}, for example) to
this operator to define a hierarchy of interpolation operators, $P_k$,
for $1 \leq k \leq k_{\text{max}}-1$.  Note that, on the finest grid,
this term corresponds to a discretized differential operator of the
form for which AMG is well-known to be effective, so no modifications
to the classical AMG heuristics are needed.  Note also that this
directly calculates the Galerkin coarse-grid operators,
$P_1^TL^T\left(D^{(\ell)}\right)^2LP_1$ and so forth, as the
coarsening algorithm proceeds.

To aid in the implementation of the multigrid V-cycle, we store
several additional pieces of data within the AMG setup phase.  First,
while the AMG setup naturally computes the hierarchy of Galerkin
coarse-grid operators, we additionally compute and store the hierarchy
of ``one-sided'' coarsening operators generated by defining $A_1 = A$,
and $A_k = A_{k-1}P_k$ as well as $M_1 = M^{(\ell)} = D^{(\ell)}L$,
and $M_k = M_{k-1}P_k$.  Note that storage of these two hierarchies of
operators gives the necessary information to compute matrix-vector
products with the full Galerkin coarsening of the normal equations
in \eqref{eq:normal} on any level of the grid hierarchy , since
$P_1^TA^TAP_1 = (AP_1)^T(AP_1)$, for example, but without the
(possibly significant) expense of forming and storing $A^TA$ (and its
Galerkin coarsenings).  Note, also, the important connection between
the one-sided coarsening hierarchies and the optimal coarse-grid
correction from the least-squares formulation, which would be to solve
\[
\min_{y_c}\left\| \left[\begin{array}{c} A \\ 
      \lambda^{(i)}D^{(\ell)}L\end{array}\right](\hat{x}+Py_c) - \left[\begin{array}{c} b \\ 0 \end{array}\right]\right\|^2,
\]
for a given approximation, $\hat{x}$, to $x^{(i,\ell)}$.  As the
interpolation and coarse-grid operators do not depend on
$\lambda^{(i)}$, these terms can be {\it formed once for each outer
iteration (corresponding to each choice of $M^{(\ell)}$) and reused
for all values of $\lambda^{(i)}$} for which a multigrid iteration is
needed. To aid in the efficient implementation of a weighted-Jacobi
style relaxation scheme (or diagonally preconditioned CG), we also compute and store
vectors corresponding to the diagonals of $A_k^TA_k$ and $M_k^TM_k$
on each level, $1 \leq k \leq k_{\text{max}}-1$.  For the latter of these, we
simply store the diagonals of the full Galerkin coarsening of
$M_1^TM_1$ on all levels.  For the former, we compute only the
diagonals of $A_k^TA_k$ by computing the sum of squares of entries in
each column of $A_k$.

While classical AMG applied to \eqref{eq:normal} would typically use a
Gauss-Seidel relaxation scheme, this is not practical here, since it
would require knowledge of $A^TA$ and its Galerkin coarse-grid
representations.  Instead, we make use of diagonally preconditioned CG
as the relaxation scheme. 
Matrix-vector products
with the normal equations on each level are computed in two stages, by
first multiplying a vector on that level, $x_k$, by $A_k$ and $M_k$ to
form $\left[\begin{smallmatrix}
A_kx_k \\ \lambda^{(i)}M_kx_k \end{smallmatrix}\right]$ and then
multiplying by their transposes to complete the matrix-vector product.
The resulting residuals can be scaled by the inverse of the diagonal
on each level by a simple vector scaling, forming
$\left(\text{diag}(A_k^TA_k)
+ \left(\lambda^{(i)}\right)^2\text{diag}(M_k^TM_k)\right)^{-1}r_k$
based on the known vectors representing $\text{diag}(A_k^TA_k)$ and
$\text{diag}(M_k^TM_k)$ and the parameter $\lambda^{(i)}$.
Within the multigrid cycle, we choose to use two sweeps of pre-relaxation and
one sweep of post-relaxation, based on studies discussed in Section \ref{ssec:relaxation_comparison}.  
The use of diagonally preconditioned CG as relaxation on each
level, instead of a stationary (or Chebyshev accelerated) weighted
Jacobi relaxation, automatically chooses an effective relaxation parameter.
This aids in the definition of a multigrid method
that is suitable for use for all values of $\lambda^{(i)}$, since no
explicit tuning of relaxation parameters as the regularization weight
changes is needed, as is shown below to be an effective approach.  The
use of Krylov acceleration within multigrid relaxation schemes is not
new \cite{REBank_CCDouglas_1985a, BReps_etal_2010a}; in preliminary
work (not reported here), experiments with fixed relaxation weights
showed that tuning was necessary with changes in both $D^{(\ell)}$ and
$\lambda^{(i)}$, so that using Krylov acceleration is an attractive
alternative.

\subsection{Improving Efficiency of Outer Iterations} \label{sec:secondimprove}
For each outermost iteration, over $M^{(\ell)}$, we must solve the
\textit{sequence} of linear systems in \eqref{eq:normal} for a set of
values, $\lambda^{(i)}$.  As a consequence of the non-stationary
relaxation procedure described above, we use Flexible GMRES
(FGMRES) \cite{YSaad_2003a} preconditioned by the multigrid V-cycle to
solve each linear system.  To do this most efficiently, we cycle
through the chosen values of $\lambda^{(i)}$ from largest to smallest,
using the solution of the next-largest regularization parameter as the
initial guess for each system.  We choose this order because the
multigrid solver is expected to be most effective for large values of
$\lambda^{(i)}$, where the regularization term dominates.  For the
largest value of $\lambda^{(i)}$, we use the solution for the chosen
$\lambda^{*,\ell-1}$ from the previous iteration as the initial guess,
with a zero initial guess used for the first outer iteration.  Note
that, since the solution is expected to vary continuously with
$\lambda^{(i)}$, this generally provides very good initial guesses, at
least as we vary the regularization parameter within each outer
iteration.  This motivates the use of a stopping tolerance that is not
relative to the initial residual; here, we require the $\ell_2$ norm
of the residual of the normal equations \eqref{eq:normal} to be less than $10^{-6}$.

\subsection{Pruning the L-curve}
Even with the efficient multigrid algorithm and improved initial
guesses, there is still a nontrivial cost to solving each linear
system in \eqref{eq:normal} for each value of $\lambda^{(i)}$.  In
order to further control the costs of the outer iterations, we now
introduce an algorithm to effectively ``trim'' the set of values
$\left\{\lambda^{(i)}\right\}$ that needs to be actively considered at
each iteration.  Recall that we solve the system for multiple values
of the parameter so that we can construct the L-curve for the
regularized problem for each choice of $M^{(\ell)}$.  As will be
observed in Section \ref{sec:numerical}, the values selected as
$\lambda^{*,\ell}$, corresponding to the ``corner'' of the L-curve,
monotonically increase with $\ell$, typically starting from the
smallest values of the set $\left\{\lambda^{(i)}\right\}$ at the first
iteration, and increasing by just one or two increments in each outer
iteration.

Rather than solving \eqref{eq:normal} for all possible values of
$\lambda^{(i)}$ at each iteration, to improve efficiency of the
algorithm, we limit ourselves to only a subset of these values.
This trimming must, clearly, be done carefully, in order to not lose
the ability to detect the corner of the L-curve from having too few
data points, or only having data points on one ``leg'' of the curve.
To do this, we consider the values $\lambda^{(i)}$ to be ordered
increasingly from smallest (denoted $\lambda^{(1)}$) to largest (with
all values assumed to be strictly positive).  For each outer
iteration, $\ell$, we keep track of the index, $i_{\ell-1}$, of the
regularization parameter chosen at iteration $\ell-1$.  Then, at iteration $\ell$, we consider only a subset
of $\left\{\lambda^{(i)}\right\}$ of size 10, ``centered'' at index
$i_{\ell-1}$.  In general (when $i_{\ell-1}$ is not too close to
either endpoint of the interval), we select those values with index
satisfying $i_{\ell-1}-2 \leq i \leq i_{\ell-1}+7$, consistent with
the expectation that $i_{\ell-1} \leq i_\ell$, but that some data for
$i \leq i_{\ell-1}$ may be needed in order to accurately pick the
corner of the L-curve.  When $i_{\ell-1}$ is close to either endpoint,
the interval is simply shifted to the smallest or largest 10 values in
$\left\{\lambda^{(i)}\right\}$.  Since it is difficult to predict good
initial values of the regularization parameter (when $\ell = 1$), and
such values vary depending on both the problem at hand and the
dimensions of $A$, we generally initialize this algorithm by using a full sweep
over the set $\lambda^{(i)}$ to find $i_1$, and trim the
values from only the second iteration onwards.  A comparison between
this strategy and one that aims to trim also at the first iteration is
given in Section \ref{ssec:trimming_numerical}, along with a
comparison to the brute-force approach of computing
all values of $x^{(i,\ell)}$.  This will show that {\it this trimming strategy can save significant
computational time without sacrificing any accuracy in the computed
solution.}  

\subsection{Stopping Criteria for Outer Iteration}\label{subsec:stopping_criterion_for_outer_iteration}

Clearly, the number of outer iterations is a large determining factor of the overall amount of work and quality of the final returned regularized solution.  
In \cite{Gazzola_etal_2020}, the heuristic for stopping the outer
iteration on $\ell$ is to track $\| L x^{(*,\ell)} \|_2$, and to terminate the outer iteration
if this
is decreasing from step $\ell-1$ to step $\ell$.  However, for certain applications
and/or high noise levels,  we found that this heuristic did not always lead to good
quality solutions.  In this subsection, we give another heuristic
that is much more robust and, thus, which we use in all of our numerical experiments.

We recall how 
the diagonal weight matrices are updated:  $D^{(\ell)} = \text{diag}(d^{(\ell)}) D^{(\ell -1)}$, where $d^{(\ell)}$ is defined from $x^{(*,\ell-1)}$ according to (\ref{eq:g}).   The weights satisfy $0 \le d^{(\ell)}_i \le 1$, so $D^{(\ell)}_{ii} \le D^{(\ell-1)}_{ii}$. 
 It is clear that $D^{(\ell)}$ will differ significantly from $D^{(\ell -1)}$ 
if we observe a previously unseen large discontinuity in $x^{(*,\ell -1)}$ (an actual edge, so $d^{(\ell)}_i = 0$ there).   
In that situation, we might expect an increase in the optimal value of the regularization parameter for the new system, and possibly a significant increase if we have succeeded to this point at capturing all the true edges (i.e., we have achieved the point at which we only need to remove noise).   
In both our experiments and in \cite{Gazzola_etal_2020}, the parameters that are selected appear to be non-decreasing as a function of $\ell$.  Indeed, the authors of the latter were able to rigorously prove that this must be the case (under mild assumptions) if the discrepancy principle, instead of the L-curve, is employed to find the parameter.

At the other extreme, note that
 if $d^{(\ell)} = {\bf 1}$ (i.e., $D^{(\ell-1)} L x^{(*,\ell-1)} = 0$), then $D^{(\ell)} = D^{(\ell-1)}$, so we would be solving the same optimization problem again, resulting in the L-curve returning the same parameter.   This extreme is unlikely to happen in practice, but we can apply similar logic:  when the changes in the weighted gradient images are small in a few consecutive iterates, we have detected all edge information and regularization is only trying to smooth noise.  Thus, the shape of the L-curve is not likely to change and the selected parameters reflect the fact that the L-curves are now very similar.

Therefore, we track the L-curve selected values $\lambda^{(*,\ell)}$ as a function of $\ell$.  When we observe that 
$\lambda^{(*,\ell-2)} = \lambda^{(*,\ell-1)} = \lambda^{(*,\ell)}$, we
stop the outer iteration at this value of $\ell$.  Requiring equal regularization parameters over three consecutive outer iterations was determined by trial-and-error.  In particular, requiring equal parameters over only two outer iterations was seen to be insufficient to determine whether there was additional edge content still to be revealed, but, after three, it seems clear there can be little further qualitative improvement in the solution.   Note that this heuristic is very cheap to employ and does not require any additional computations.

\subsection{Resulting method}
Overall, we obtain the method provided
in \Cref{alg:improved_regularization_using_multigrid}, assuming that
$m\geq 10$ values are provided in the set of regularization
parameters, $\{\lambda^{(i)}\}_{i=1}^m$.
\begin{algorithm}[h]
  \caption{Improved regularization using FGMRES-accelerated multigrid}
  \label{alg:improved_regularization_using_multigrid}
  \begin{algorithmic}
  \FOR{$\ell = 1$ to maxits}
  	\IF{$\ell == 1$}
		\STATE{$lam\_range \leftarrow \{\lambda^{(1)},\dots,\lambda^{(m)}\}$}
	\ELSE
		\STATE{$lam\_range \leftarrow \{\lambda^{(\max\{\min\{i_{\ell-1}-2,m-9\},1\})},\dots,\lambda^{(\min\{\max\{i_{\ell-1}+7,10\},m\})}\}$}
	\ENDIF
	\STATE{setup AMG}
	\FOR{$\lambda^{(i)}$ in $lam\_range$}
		\STATE{solve $A^TA x_i + (\lambda^{(i)})^2 L^T (D^{(\ell)})^2 L x^{(i,\ell)} = A^T b$ using FGMRES-accelerated AMG}
	\ENDFOR
	\STATE{select $\lambda^{(*,\ell)}$ using L-curve, $i_\ell \leftarrow \text{index of chosen $\lambda$}$}
  	\IF{$\ell>2$ and $\lambda^{(*,\ell-2)} == \lambda^{(*,\ell-1)} == \lambda^{(*,\ell)}$}
		\STATE{return $x^{(i_{\ell},\ell)}$}
        \ELSE
                \STATE{compute $D^{(\ell+1)}$}
	\ENDIF
  \ENDFOR
  \end{algorithmic}
\end{algorithm}

\section{Numerical Results}  \label{sec:numerical}
The described method has been implemented and tes\-ted using MATLAB R2022a and the toolboxes AIR Tools II \cite{art:HANS18} and IR Tools \cite{art:GAZZ19}. We stop the outer iteration when
\[
\lambda^{(*,\ell-2)} = \lambda^{(*,\ell-1)} = \lambda^{(*,\ell)},
\]
as described above in \Cref{subsec:stopping_criterion_for_outer_iteration}. The multigrid preconditioned
FGMRES iteration was used (without restarting), with a stopping
criterion of a reduction of the $\ell^2$ norm of the residual of the normal equations to below $10^{-6}$. As possible values for the regularization parameter, $\lambda$, in the outer iteration, we used 30 values equally distributed on a logarithmic scale between $10^{-3}$ and $10^{2}$ using MATLAB's \texttt{logspace(2,-3,30)} command.

\subsection{Usage of the solution of the previous inner iteration}
\begin{figure}[htbp]
\begin{center}
\includegraphics{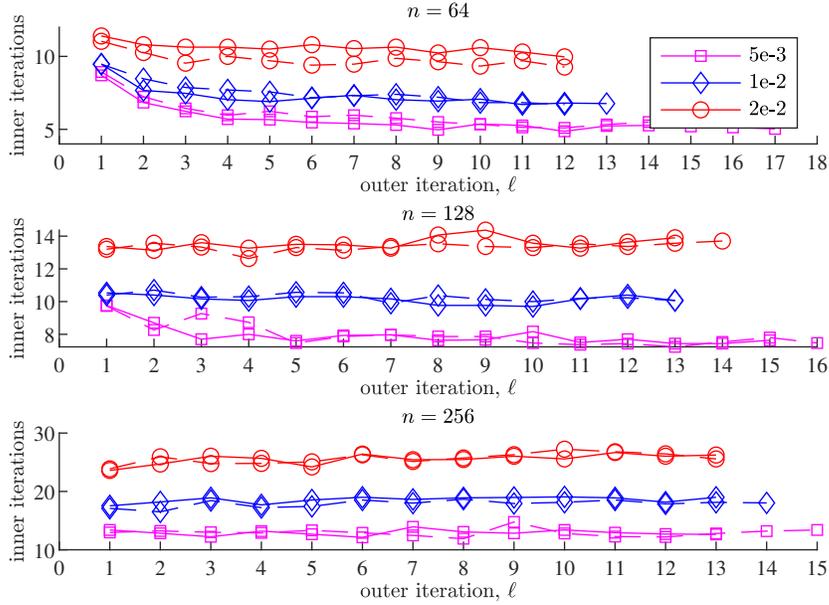}%
\caption{Average number of inner iterations needed for the different lambdas in each outer iteration of the method applied to the Shepp Logan phantom with zero initial guess (solid lines) or by using the solution from the previous iteration (dashed line).}
\label{fig:inner_iteration_zero_initial_guess}
\end{center}
\end{figure}

Before evaluating the proposed method in detail, we want to highlight the effect of taking the initial guess for each system in the inner iteration as the solution of the previous inner iteration.  In Fig.~\ref{fig:inner_iteration_zero_initial_guess}, we compare two strategies: either using a zero initial guess for each linear system to be solved in the inner iteration process or, starting from the second inner iteration in each outer iteration, using the solution of the previous inner iteration as the initial guess. We consider a tomography example using the Shepp Logan phantom of different sizes with different amounts of added noise, created using IR Tools' \cite{art:GAZZ19} \texttt{PRtomo} and \texttt{PRnoise} routines. The resulting average number of inner iterations per outer iteration is plotted in Fig.~\ref{fig:inner_iteration_zero_initial_guess}. While the effect is small, it demonstrates that using the available information can reduce the number of inner iterations necessary to reach convergence, especially for larger error levels, where the inner iteration counts were highest.


\subsection{Effect of pre- and post-relaxation}
\label{ssec:relaxation_comparison}
\begin{figure}[htbp]
\begin{center}
\includegraphics{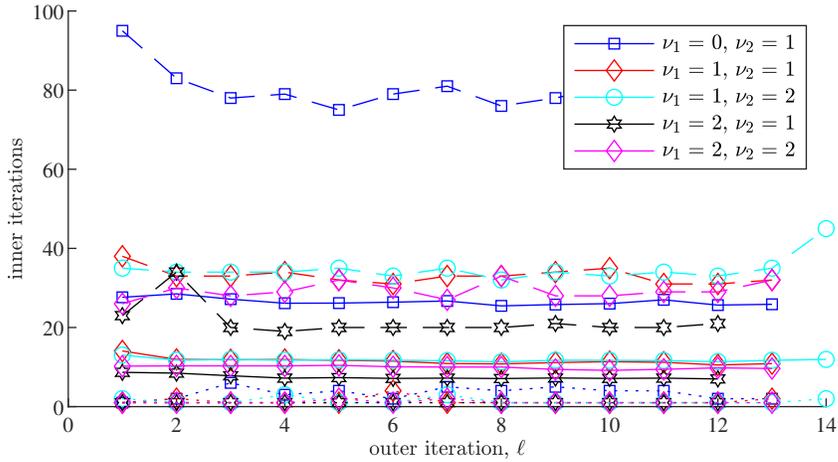}
\caption{Number of inner iterations needed in each outer iteration of the method applied to the Shepp Logan phantom of size $128 \times 128$ with 1\% noise for different multigrid cycles. Solid lines represent the average number of iterations needed for the different values of lambda, dashed lines the maximum number of iterations, and dotted lines the minimum.}
\label{fig:its_cycling}
\end{center}
\end{figure}

\begin{figure}[htbp]
\begin{center}
\includegraphics{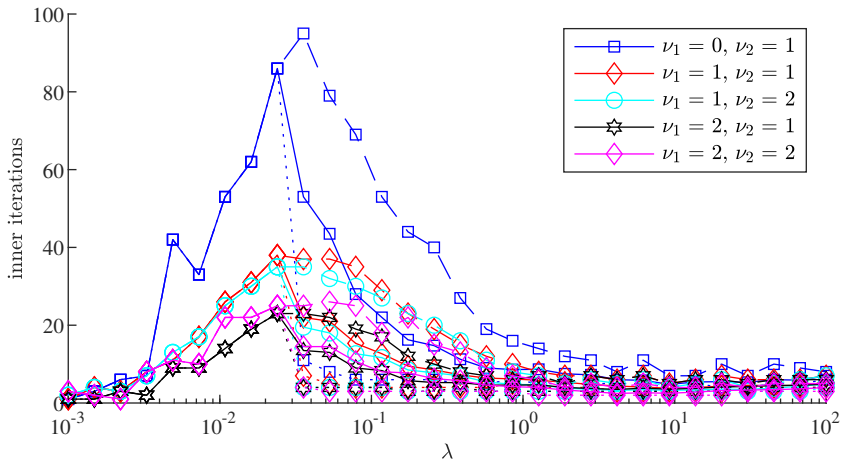}
\caption{Number of inner iterations needed for each lambda when the method is applied to the Shepp Logan phantom of size $128 \times 128$ with 1\% noise for different multigrid cycles. Solid lines represent the average number of iterations needed for the different values of lambda, dashed lines the maximum number of iterations, and dotted lines the minimum.}
\label{fig:its_cycling_lambda}
\end{center}
\end{figure}

An important parameter that has to be chosen in multigrid methods is the number of pre- and post-relaxation steps. To be able to do this, we study the dependence of the convergence rate on the number of pre- and post-relaxation steps. For that purpose, we study the same tomography example as before, limiting ourselves to size $128 \times 128$ with 1\% added noise. We tested V(0,1)-, V(1,1)-, V(1,2)-, V(2,1)- and V(2,2)-cycles for this problem. For each outer iteration, we recorded the minimal, maximal and average number of inner iterations needed for the different values of $\lambda$ that were evaluated in the respective outer iteration. The resulting maximal, minimal and average numbers of iterations are depicted in Fig.~\ref{fig:its_cycling}, while the average number of iterations needed in each outer iteration can also be found in Table~\ref{tab:its_cycling}. 
Additionally, in Fig.~\ref{fig:its_cycling_lambda} the numbers of iterations for the different lambdas is shown.
As can be seen from these results, the difference in performance between V(2,1)-cycles and V(2,2)-cycles is small, and the former even perform better. 
On the other hand, the improvement given by using V(2,1)-cycles over V(1,1)- or V(1,2)-cycles is clear.  
Therefore, we use V(2,1) cycles in the rest of our numerical experiments.
\begin{table}[htp]
\caption{Average number of inner iterations needed in each outer iteration of the method applied to the Shepp Logan phantom of size $128 \times 128$ with 1\% noise for different multigrid cycles.}
\begin{center}
\begin{tabular}{|l|r|r|r|r|r|r|r|}
\hline
& \multicolumn{1}{c|}{1} & \multicolumn{1}{c|}{2} & \multicolumn{1}{c|}{3} & \multicolumn{1}{c|}{4} & \multicolumn{1}{c|}{5} & \multicolumn{1}{c|}{6} & \multicolumn{1}{c|}{7} \\
\hline
V(0,1) & 27.63 & 28.50 & 27.17 & 26.13 & 26.17 & 26.40 & 26.70 \\
V(1,1) & 14.07 & 11.93 & 11.83 & 11.90 & 11.70 & 11.50 & 10.90 \\
V(1,2) &13.03 & 11.73 & 11.93 & 11.73 & 11.83 & 11.73 & 11.60 \\
V(2,1) & 8.67 & 8.47 & 7.80 & 7.20 & 7.33 & 7.13 & 7.20 \\
V(2,2) & 8.67 & 8.47 & 7.80 & 7.20 & 7.33 & 7.13 & 7.20 \\
\hline
\hline
& \multicolumn{1}{c|}{8} & \multicolumn{1}{c|}{9} & \multicolumn{1}{c|}{10} & \multicolumn{1}{c|}{11} & \multicolumn{1}{c|}{12} & \multicolumn{1}{c|}{13} & \multicolumn{1}{c|}{14} \\
\hline
V(0,1) & 25.47 & 25.80 & 26.03 & 27.00 & 25.67 & 25.87 & -- \\
V(1,1) & 10.80 & 11.13 & 11.37 & 11.20 & 10.53 & 10.87 & -- \\
V(1,2) & 11.30 & 11.73 & 11.70 & 11.63 & 11.37 & 11.70 & 11.93 \\
V(2,1) & 7.10 & 7.23 & 7.13 & 7.20 & 7.00 & -- & -- \\
V(2,2) & 7.10 & 7.23 & 7.13 & 7.20 & 7.00 & -- & -- \\
\hline
\end{tabular}
\end{center}
\label{tab:its_cycling}
\end{table}%

\subsection{Effect of trimming the L-curve}
\label{ssec:trimming_numerical}
\begin{figure}[htbp]
\begin{center}
\includegraphics{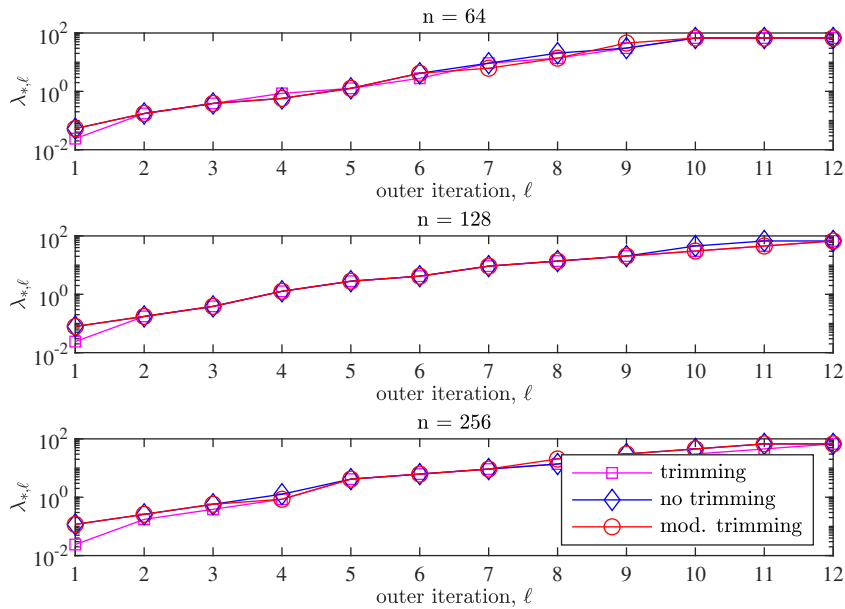}%
\caption{Effect of trimming the L-curve for a full angle computer tomography using the Shepp Logan phantom with 1\% noise. Depicted are the regularization parameters chosen by the algorithm. The modified trimming uses trimming only after the first iteration.}%
\label{fig:L-Curve_trimming}%
\end{center}
\end{figure}

Using the V(2,1)-cycle and 1\% noise, we next consider the effect of the proposed trimming of the L-curve for a computer tomography example using the Shepp Logan phantom of sizes $64 \times 64$, $128 \times 128$ and $256 \times 256$. We compare the L-curve algorithm from \cite{HansenOleary93} applied to the whole range of values of $\lambda$ with two variants that apply it to a smaller range of values. The first variant uses only the smallest 10 values of $\lambda$ when computing the solution to the first outer iteration and, after that, takes the range described above with the seven larger values of $\lambda$ than the optimal value chosen at the previous outer iteration, that value itself, and the two immediately smaller values. The second variant uses the whole range of $\lambda$ values at the first outer iteration but, from the second iteration onwards, follows the same strategy as the first variant. The regularization parameters chosen are depicted in Fig.~\ref{fig:L-Curve_trimming}. Here, we see that, for these examples, the first variant tends to choose a too small regularization parameter in the first iteration. The second variant, in contrast, chooses almost exactly the same regularization parameters and shows the same algorithmic behavior as the variant without trimming, but at a little more than one third of the computational cost.
\subsection{Full angle computer tomography}
\begin{figure}[htbp]
\begin{center}
\includegraphics{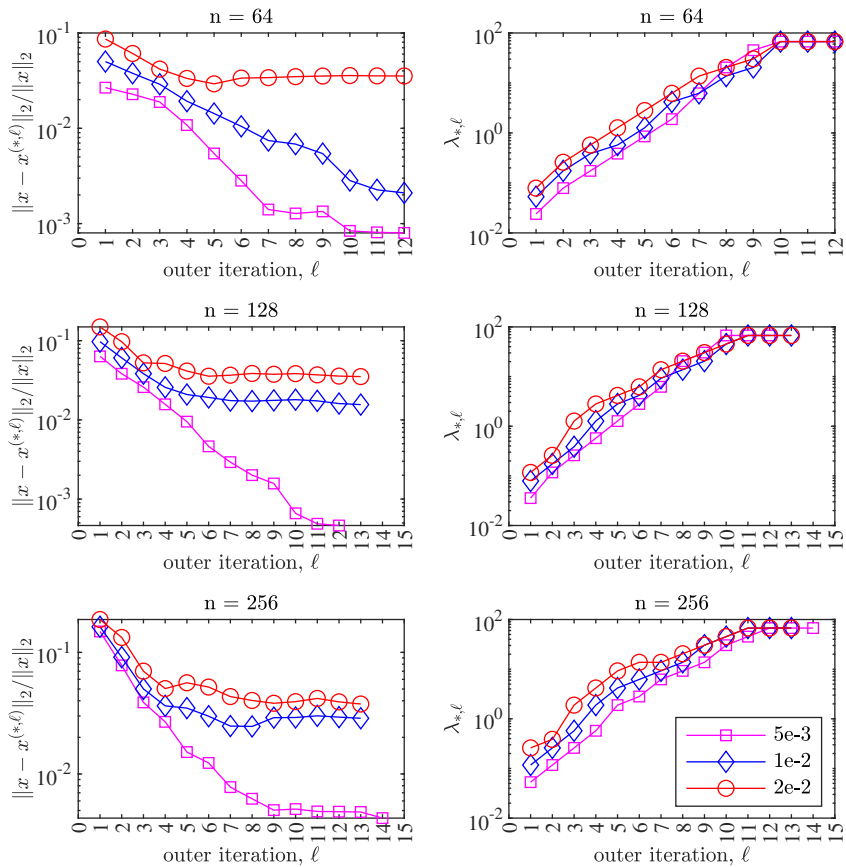}%
\caption{Errors (left) and regularization parameters chosen (right) for full angle computer tomography of the Shepp Logan phantom for different image sizes (by row) and noise levels.}%
\label{fig:full_angle_errs_and_reg_params}%
\end{center}
\end{figure}

\begin{figure}[htbp]
\begin{center}
\includegraphics{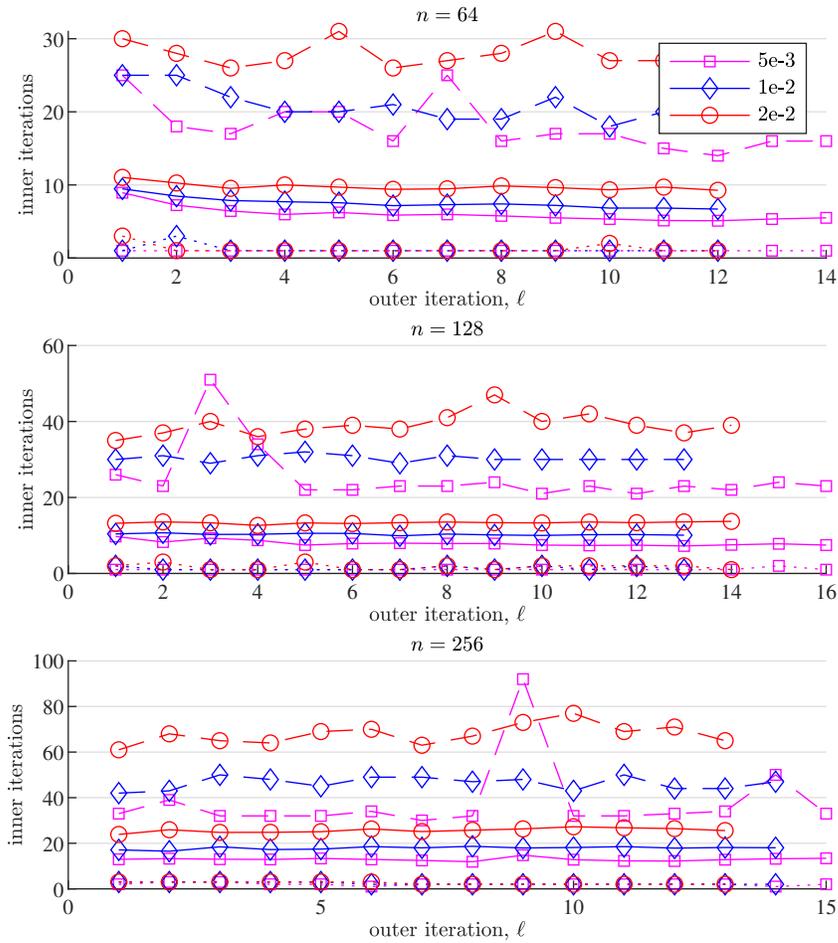}%
\caption{Inner iterations for full angle computer tomography of the Shepp Logan phantom.  Solid lines represent the average number of iterations needed for the different values of lambda, dashed lines the maximum number of iterations, and dotted lines the minimum.}
\label{fig:full_angle_inner_its}
\end{center}
\end{figure}

\begin{figure}[htbp]
\begin{center}
\includegraphics{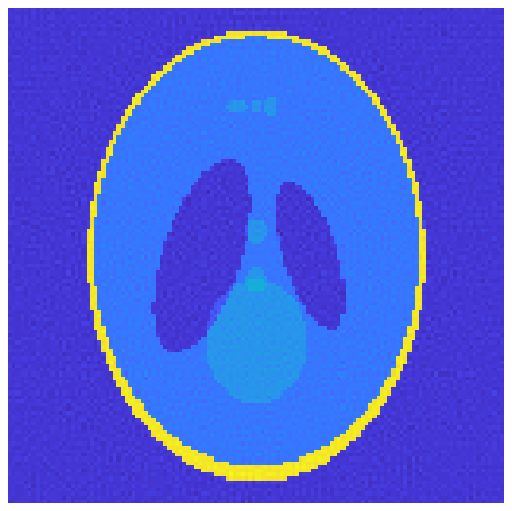}\includegraphics{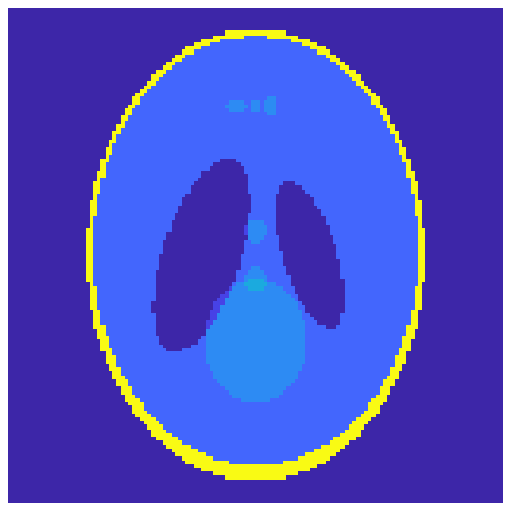}
\caption{Initial and final image obtained for the Shepp Logan computer tomography example for size $128 \times 128$ and 0.5\% noise.}
\label{fig:full_angle_initial_and_final}
\end{center}
\end{figure}

For the computer tomography example using the Shepp Logan phantom of sizes $64 \times 64$, $128 \times 128$ and $256 \times 256$, we now compare the performance of the algorithm for different noise levels. We used the modified trimming approach and a V(2,1)-cycle for preconditioning FGMRES for 0.5\% noise, 1\% noise, and 2\% noise. The resulting errors and the chosen regularization parameters can be found in Fig.~\ref{fig:full_angle_errs_and_reg_params}.
The minimal, maximal and average inner iterations needed can be found in Fig~\ref{fig:full_angle_inner_its}. 
The initial and final reconstructions are in Fig.~\ref{fig:full_angle_initial_and_final}.

In general, results here are as we expect them to be. Fig.~\ref{fig:full_angle_errs_and_reg_params} 
shows that errors in $x^{(*,\ell)}$ generally decay with each outer iteration, and that the outer iteration stopping criterion prevents significant issues with oversmoothing of the resulting image (as shown, for example, in Fig.~\ref{fig:full_angle_initial_and_final}).  Fig~\ref{fig:full_angle_inner_its} shows that inner iteration counts are relatively consistent over the course of the outer iterations, although they are clearly sensitive to noise level and problem size.
\subsection{Limited angle computer tomography}
\begin{figure}[htbp]
\begin{center}
\includegraphics{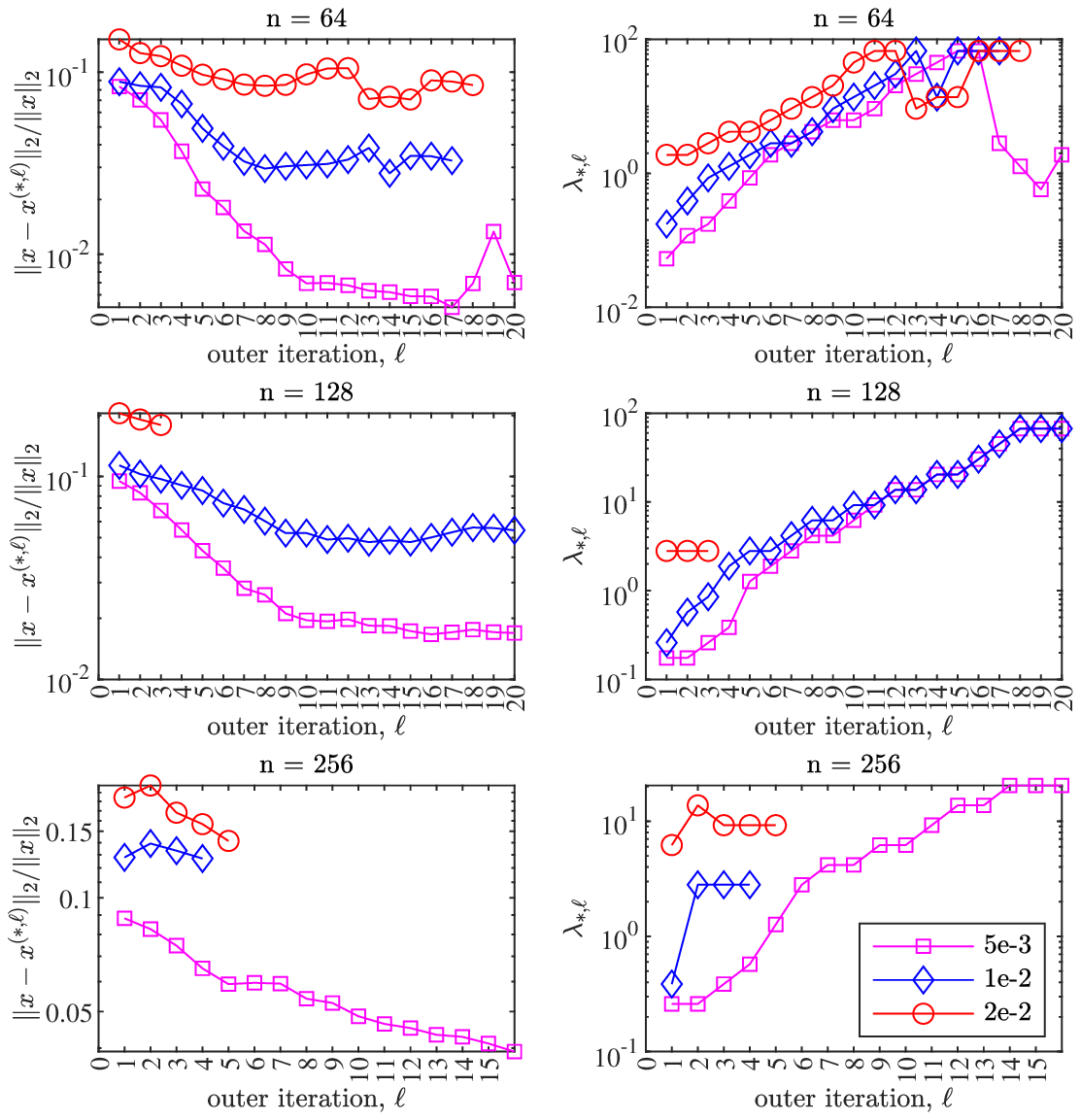}%
\caption{Errors (left) and regularization parameters chosen (right) for the 130-degree limited angle computer tomography of the grains example for different image sizes (by row) and noise levels.}
\label{fig:limited_angle_grains_errs_and_reg_params}
\end{center}
\end{figure}

\begin{figure}[htbp]
\begin{center}
\includegraphics{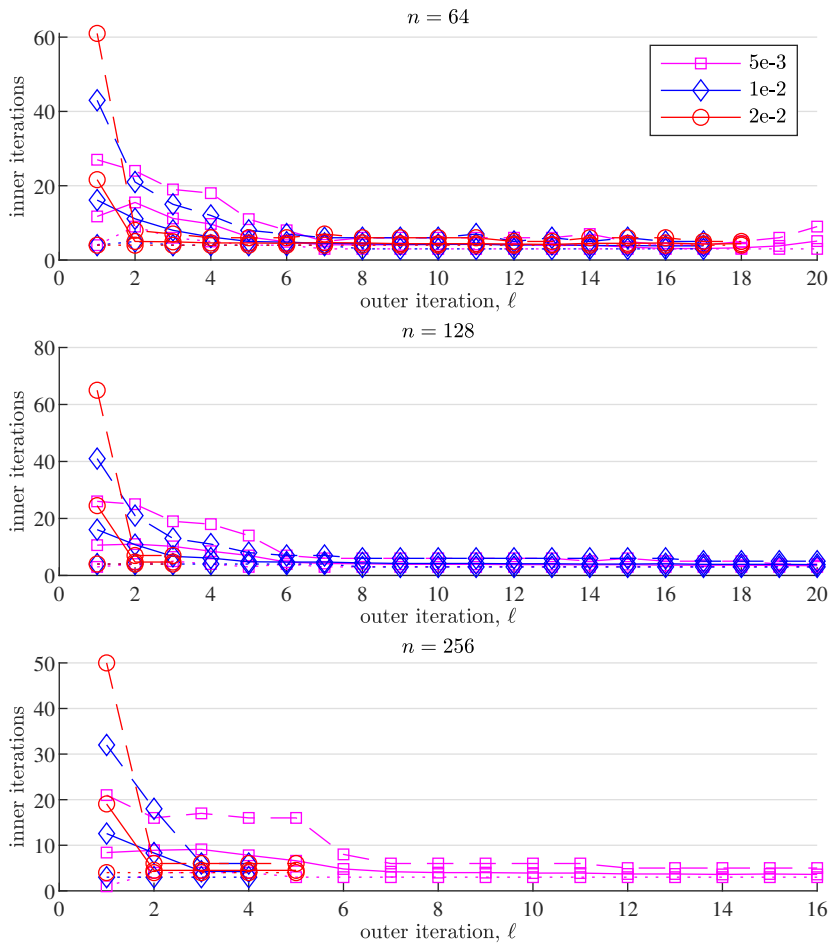}%
\caption{Inner iterations for limited angle computer tomography of the grains example.  Solid lines represent the average number of iterations needed for the different values of lambda, dashed lines the maximum number of iterations, and dotted lines the minimum.}
\label{fig:limited_angle_inner_its}
\end{center}
\end{figure}

\begin{figure}[htbp]
\begin{center}
\includegraphics{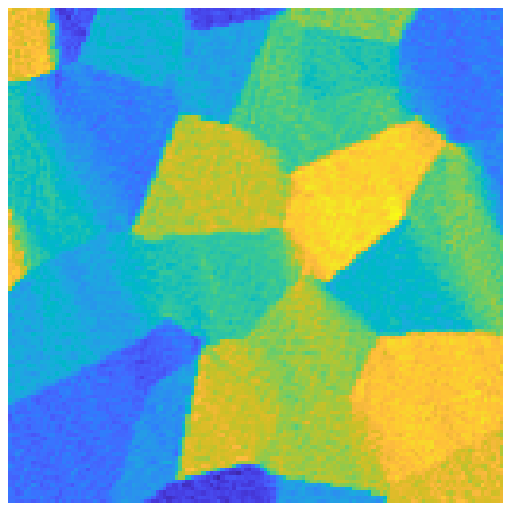}\includegraphics{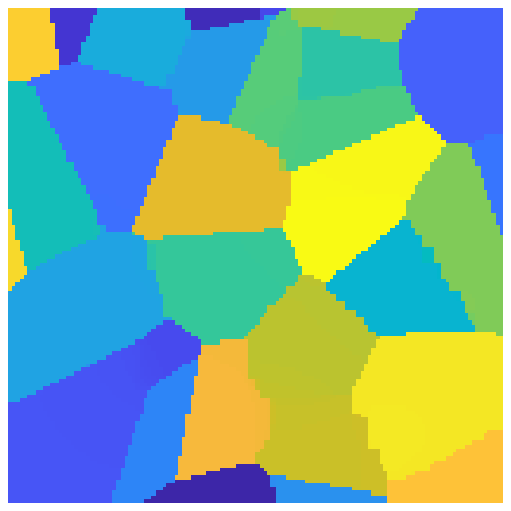}
\caption{Initial and final image obtained for the grains limited angle computer tomography example for size $128 \times 128$ and 0.5\% noise.}
\label{fig:limited_angle_grains_initial_and_final}
\end{center}
\end{figure}

We now consider a similar example for \textit{limited angle} computer tomography, using the grains example from IR tools \cite{art:GAZZ19}. We use the command \texttt{PRset('angles', 0:2:130, 'phantomImage', 'grains')} to set up the problem and choose sizes $64 \times 64$, $128 \times 128$, and $256 \times 256$. As before, we compare the performance using the modified trimming and a V(2,1)-cycle for preconditioning FGMRES with 0.5\% noise, 1\% noise, and 2\% noise. The minimal, maximal and average inner iterations needed can be found in Fig~\ref{fig:limited_angle_inner_its}. 
Fig.~\ref{fig:limited_angle_grains_errs_and_reg_params} 
shows the resulting errors and the chosen regularization parameters. The initial and final reconstructions for the $128\times 128$ image with 0.5\% noise are shown in Fig.~\ref{fig:limited_angle_grains_initial_and_final}.  As above, we find these results are mostly as expected, showing consistent reduction in error norm with outer iterations.

\subsection{Image deblurring}
\begin{figure}[htbp]
\begin{center}
\includegraphics{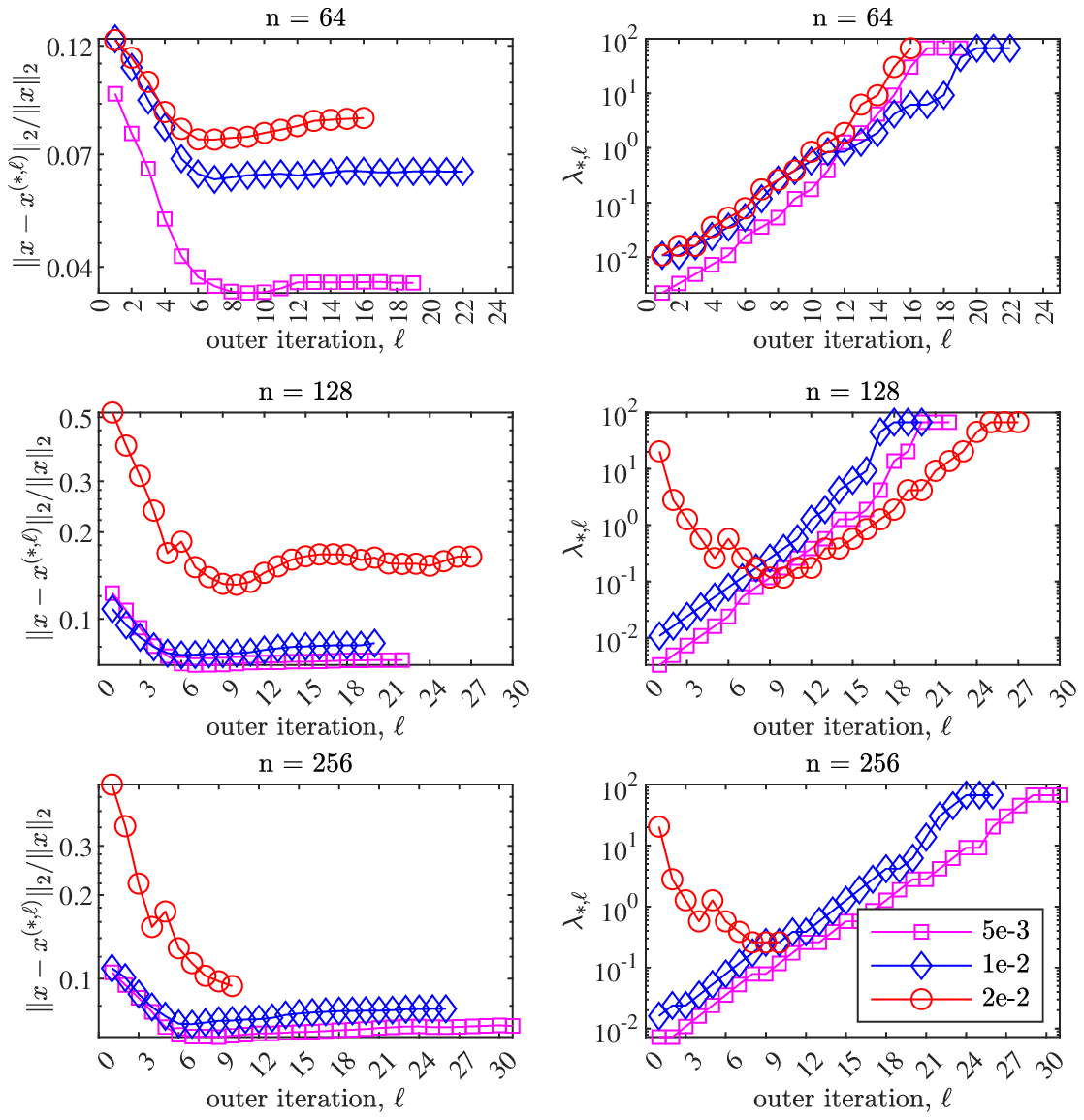}%
\caption{Errors (left) and regularization parameters chosen (right) for the blurred grains example for different image sizes (by row) and noise levels.}
\label{fig:grains_blurring_errs_and_reg_params}
\end{center}
\end{figure}

\begin{figure}[htbp]
\begin{center}
\includegraphics{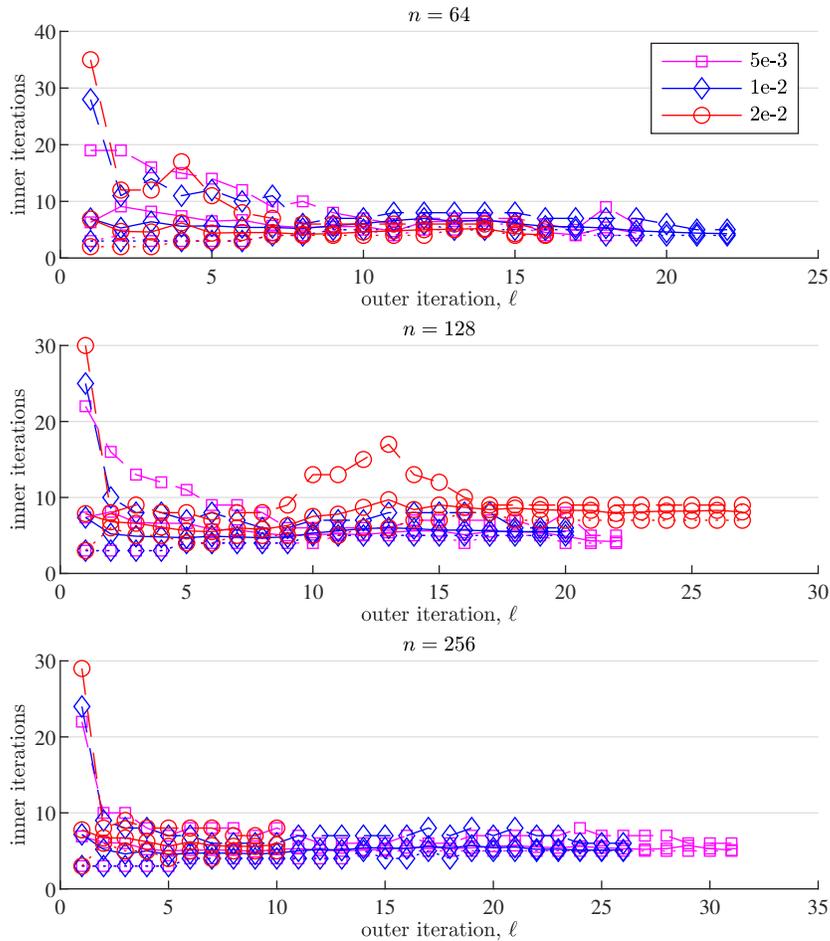}%
\caption{Inner iterations for the blurred grains example.  Solid lines represent the average number of iterations needed for the different values of lambda, dashed lines the maximum number of iterations, and dotted lines the minimum.}
\label{fig:grains_blurring_inner_its}
\end{center}
\end{figure}

\begin{figure}[htbp]
\begin{center}
\includegraphics{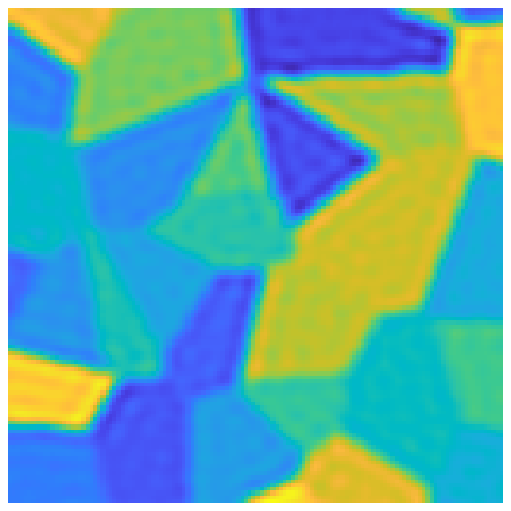}\includegraphics{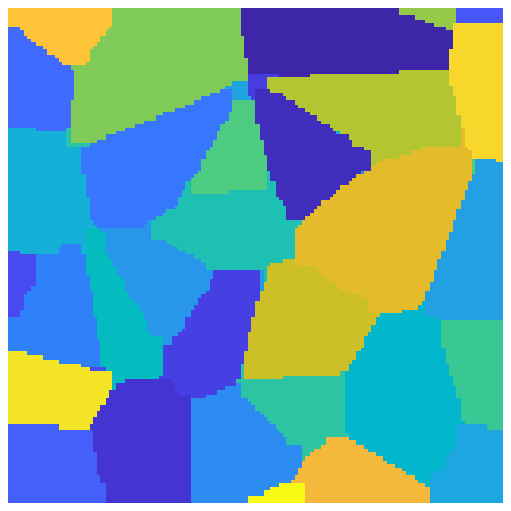}
\caption{Initial and final image obtained for the blurred grains example for size $128 \times 128$ and 0.5\% noise.}
\label{fig:grains_blurring_initial_and_final}
\end{center}
\end{figure}

To demonstrate that the methodology also works for other inverse problems, we now blur the grains example images of sizes $64 \times 64$, $128 \times 128$ and $256 \times 256$, using a discretized integral kernel for the $m \times n$ image represented by the Toeplitz matrix created using:
\begin{verbatim}

band1 = 8; sigma1=1.5*n/64; 
band2 = 7; sigma2=1.25*n/64;
z1 = [exp(-((0:band1-1).^2)/(2*sigma1^2)),zeros(1,n-band1)];
z1 = z1./(2*sum(z1)-1); z1=z1';
z2 = [exp(-((0:band2-1).^2)/(2*sigma2^2)),zeros(1,m-band2)];
z2 = z2./(2*sum(z2)-1); z2=z2';
A1 = toeplitz(z1); A1 = sparse(A1); 
A2 = toeplitz(z2); A2 = sparse(A2);
A = kron(A1,A2);

\end{verbatim}
As before, 0.5\%, 1\%, and 2\% noise are added to the blurred image and the method is applied using the same settings. The behavior of the algorithm is similar to the other cases and errors as well as regularization parameters are to be found in Fig.~\ref{fig:grains_blurring_errs_and_reg_params}.
The minimal, maximal and average inner iterations needed can be found in Fig~\ref{fig:grains_blurring_inner_its}. 
The initial and final reconstruction for the $128\times 128$ image with 0.5\% noise are shown in Fig.~\ref{fig:grains_blurring_initial_and_final}.

\section{Conclusions}  \label{sec:conclusions}
In this paper, we develop a multi-pronged approach to improve the computational efficiency of solving a sequence of  
regularized systems of the form~\eqref{eq:LSreg} when $M^{(\ell)}$ takes the form of a diagonally weighted gradient operator changing with each outer iteration, $\ell$.  
  When $\lambda$ is large, we use the similarity of the corresponding regularized least squares problems to those involving diffusion operators to develop an AMG-style preconditioner.  The bulk of the work in precomputing quantities for the preconditioner is done once for each $\ell$ and reused across each $\lambda$, making it very efficient to implement in practice. To account for the rectangular structure---and the fact that the discrete operator, $A$, will usually possess many more non-zeros per row than the discrete differential operators to which multigrid is usually applied---the preconditioner is built for the normal equations and the resulting prolongation operators are applied from the right to $A$ and $L$ to obtain a smaller rectangular system. This allows us to  set up suitable coarse-level problems while keeping the computational complexity relatively low.

Before finding the optimal regularized solution for each $\ell$, we must determine a suitable regularization parameter.  To further reduce the number of solves that are necessary for a particular $\ell$, we also develop a pruning strategy for the L-curve.  Our approach preserves the necessary structure (i.e.\ corner) of the L-curves, allowing for a good parameter to be selected at each outer iteration, even as the shape of the L-curve changes (due to changes in the regularization operator) as a function of $\ell$.  We also incorporate initial guesses, in the form of solutions to neighboring regularized problems, though we find that the required number of inner iterations is relatively insensitive to this choice.  Finally, we provide a cheap heuristic for stopping the outer iterations, to avoid unnecessary computation, and demonstrate its effectiveness by employing it in all our numerical experiments.  
We demonstrate our strategies on problems in both CT image reconstruction and image restoration of varying size.  The results reveal that the pruning strategy still allows for the 'correct' selection of the regularization parameter, irrespective of noise level and problem dimension. In all cases, the preconditioner effectively reduces the number of iterations needed over naive approaches, varying only a little over a wide range of values of $\lambda$.

The case for our innovations relative to the naive approach is clear.  A detailed numerical study against the hybrid solver approach for various classes of problems is the subject of future research, requiring optimized implementations of both approaches.   We note that, in this work, we use the developed AMG preconditioner for all values of $\lambda$, even though the greatest benefit is obtained for relatively large values.  For small values of $\lambda$, recent work~\cite{KilmerSturlerOconnell20} has shown that subspace recycling can be used very effectively.  It will be very interesting to combine subspace recycling with preconditioning and pruning techniques proposed herein, in order to achieve maximum efficiency.  
We also note that we have considered only one type of edge-preserving regularizer to motivate our work here.  However, the methods that we propose should be readily adaptable for different forms of edge-preserving regularization.  In future work, we will consider other extensions, for example, in dynamic imaging, where the regularization depends on gradient information in both space and time.

\bibliographystyle{siam} \bibliography{regularized}

\end{document}